\documentclass[10pt]{article}
\usepackage{cite}
\usepackage{mathrsfs}
\usepackage{amsfonts}
\usepackage{amsmath}
\usepackage{amsfonts,amssymb,color}
\usepackage{dsfont}
\usepackage{curves}
\usepackage{mathrsfs}
\usepackage{pifont}
\usepackage{amssymb}
\allowdisplaybreaks

\numberwithin{equation}{section}

\date{}

\def\BigRoman{\uppercase\expandafter{\romannumeral\number\count 255 }}
\def\Romannumeral{\afterassignment\BigRoman\count255=}

\begin{document}
\title{A spectral condition for a graph to have strong parity factors
}
\author{\small  Sizhong Zhou$^{1}$\footnote{Corresponding
author. E-mail address: zsz\_cumt@163.com (S. Zhou)}, Tao Zhang$^{2}$, Qiuxiang Bian$^{1}$\\
\small 1. School of Science, Jiangsu University of Science and Technology,\\
\small Zhenjiang, Jiangsu 212100, China\\
\small 2.  School of Economics and management, Jiangsu University of Science and Technology,\\
\small Zhenjiang, Jiangsu 212100, China
}

\maketitle
\begin{abstract}
\noindent A graph $G$ has the strong parity property if for every subset $X\subseteq V(G)$ with $|X|$ even, $G$ has a spanning subgraph
$F$ satisfying $\delta(F)\geq1$, $d_F(u)\equiv1$ (mod 2) for any $u\in X$, and $d_F(v)\equiv0$ (mod 2) for any $v\in V(G)\setminus X$. In this
paper, we give a spectral radius condition to guarantee that a connected graph has the strong parity property.
\\
\begin{flushleft}

(2020) Mathematics Subject Classification: 05C50, 05C70
\end{flushleft}
\end{abstract}

\section{Introduction}

In this paper, we consider finite simple graphs without loops and multiple edges. Let $G$ be a graph with vertex set $V(G)$ and edge
set $E(G)$. The number of vertices in $G$ is called its \emph{order} and denoted by $|V(G)|$. The \emph{degree} $d_G(v)$ of a vertex $v$ in a graph $G$
is the number of edges of $G$ incident with the vertex $v$. The minimum degree of the vertices of $G$ is denoted by $\delta(G)$. For any subset $S$ of $V(G)$, the subgraph of $G$ induced by $S$ is denoted by $G[S]$. We write $G-S=G[V(G)\setminus S]$
and denote by $c(G)$ the number of components in $G$. Let $K_n$ denote the complete graph of order $n$.

Let $G_1$ and $G_2$ be two vertex disjoint graphs. We denote by $G_1\cup G_2$ the \emph{disjoint union} of $G_1$ and $G_2$. For any positive integer
$t$, we denote by $tG$ the disjoint union of $t$ copies of $G$. The \emph{join} $G_1\vee G_2$ is the graph derived from $G_1\cup G_2$ by adding all
possible edges between $V(G_1)$ and $V(G_2)$.

Given a graph $G$ of order $n$, the \emph{adjacency matrix} of $G$ is defined as the matrix $A(G)=(a_{ij})_{n\times n}$ with $a_{ij}=1$ if $v_i$ and
$v_j$ are adjacent, and $a_{ij}=0$ otherwise. The largest eigenvalue of $A(G)$ is denoted by $\lambda_1(G)$. Note that the \emph{adjacency spectral
radius} (or \emph{spectral radius}, for short) of $G$ is equal to $\lambda_1(G)$, written as $\rho(G)$.

Let $g,f$ be two functions from $V(G)$ to the nonnegative integers with $g(x)\leq f(x)$ for every $x\in V(G)$. A \emph{$(g,f)$-factor} of $G$ is a
spanning subgraph $F_1$ of $G$ satisfying $g(x)\leq d_{F_1}(x)\leq f(x)$ for every vertex $x$ of $G$. A \emph{$(g,f)$-parity factor} of $G$ is a
spanning subgraph $F_2$ satisfying $d_{F_2}(v)\equiv g(v)\equiv f(v)$ (mod 2) and $g(v)\leq d_{F_2}(v)\leq f(v)$ for every $v\in V(G)$. If $g(v)=1$
and $f(v)=b$ for each $v\in V(G)$, then a $(g,f)$-parity factor is a \emph{$(1,b)$-odd factor}. Note that a $(1,b)$-odd factor is an extension of a
perfect matching (i.e., 1-factor). For any real function $\varphi$ defined on $V(G)$ and any subset $S\subseteq V(G)$, write
$\varphi(S)=\sum\limits_{v\in S}{\varphi(v)}$.

The study on factors and parity factors in graphs attracted much attention. Nishimura \cite{N}, Niessen, and Randerath
\cite{NR}, Enomoto et al. \cite{EJKS} provided some sufficient conditions for the existence of $k$-factors in graphs. Many
researchers \cite{ZZS,ZSL1,WZhi} verified some results related to binding number, degree condition, neighborhood etc., for a graph
having a $(1,2)$-factor. Many efforts have been devoted to finding sufficient conditions for graphs to have $(a,b)$-factors by utilizing various graphic
parameters such as independence number \cite{KL} and Fan-type condition \cite{Mf}.

Amahashi \cite{A} presented a characterization for a graph with a $(1,b)$-odd factor. Kim et al. \cite{KOPR} obtained an eigenvalue condition for graphs to have $(1,b)$-odd factors. Zhou, and Liu
\cite{ZL} established a lower bound on the spectral radius (resp. the size) of a graph $G$ to ensure that $G$ has a $(1,b)$-odd factor. Liu, and Lu
\cite{LL} showed a degree condition for the existence of an $(a,b)$-parity factor in a graph. Yang et al. \cite{YZLL} provided an independence
number and connectivity condition for a graph to contain an $(a,b)$-parity factor.

A graph $G$ has the strong parity property if for every subset $X\subseteq V(G)$ with $|X|$ even, $G$ has a spanning subgraph $F$ satisfying
$\delta(F)\geq1$, $d_F(u)\equiv1$ (mod 2) for any $u\in X$, and $d_F(v)\equiv0$ (mod 2) for any $v\in V(G)\setminus X$. Bujt\'as, Jendrol', and Tuza
\cite{BJT} introduced the definition of strong parity factor and provided some sufficient conditions for the existence of strong parity factors in
graphs. Lu, Yang, and Zhang \cite{LYZ} obtained a necessary and sufficient condition for a graph to have the strong parity property. Hasanvand
\cite{H} investigated the existence of strong parity factors in terms of edge-connected graphs. Zhou, and Zhang \cite{ZZ} put forward two sufficient
conditions for a graph to have the strong parity property.

O \cite{O} gave a spectral radius condition for a graph to contain a perfect matching. Li, and Miao \cite{LM} established a relationship between the
spectral radius and $(1,2)$-factors of graphs.

Motivated by \cite{O,LYZ} directly, it is natural and interesting to find a sufficient spectral condition to guarantee that a graph has the strong
parity property. In this paper, we establish a connection between the spectral radius and the strong parity property of a connected graph. Let
$H(n,\delta,s)=K_s\vee(K_{n-(\delta-1)s-1}\cup((\delta-2)s+1)K_1)$ and $H(n,\delta)=K_{\delta}\vee(K_{n-\delta(\delta-1)-1}\cup(\delta(\delta-2)+1)K_1)$,
where $s$ is a positive integer.

\medskip

\noindent{\textbf{Theorem 1.1.}} Let $G$ be a connected graph of order $n$ with minimum degree $\delta\geq3$. If $n\geq2\delta^2$ and
$$
\rho(G)\geq\rho(H(n,\delta)),
$$
then $G$ has the strong parity property unless $G=H(n,\delta)$.

\medskip

The proof technique in the paper for Theorem 1.1 follow the ideas of Li, and Miao \cite{LM}, Zhou, Sun, and Liu \cite{ZSL1}, Shen \cite{Sh}. The proof
employs standard techniques in spectral graph theory, such as the use of eigenvalues and eigenvectors, equitable partitions, and the application of
known results. But some new techniques (such as properties of similar matrices, properties of functions and derivative functions, and use of the well-known
Cauchy Interlacing Theorem, and so on) in the proof of Theorem 1.1 are used. The paper extends the existing body of knowledge on strong parity factors by
focusing on the spectral radius. This novel approach can enrich the spectral graph theory and its applications in studying the structural properties of
graphs.

In what follows, we draw a diagram for showing the spectral radius of the graph ``$H(n,\delta)$", for small $n$ and $\delta$, to show how much it is
far from the spectral radius of the complete graph.

\begin{tabular}{|c|c|c|}
\hline
$n,\delta$&$\rho(H(n,\delta))$ is probably close to &$\rho(K_n)$\\
\hline
$n=18,\delta=3$&13&17\\
\hline
$n=32,\delta=4$&22&31\\
\hline
$n=50,\delta=5$&33&49\\
\hline
$n\geq2\delta^{2},\delta\geq3$&$n-\delta(\delta-2)$&$n-1$\\
\hline
\end{tabular}

\medskip

A subgraph $F\subseteq G$ is called a \emph{positive factor} if $F$ is a spanning subgraph and has minimum degree $\delta(F)\geq1$. Bujt\'as, Jendrol',
and Tuza \cite{BJT} showed a sufficient condition for a graph to have the strong parity property.

\medskip

\noindent{\textbf{Theorem 1.2}} (Bujt\'as, Jendrol', and Tuza \cite{BJT}). Let $G$ be a connected graph with minimum degree $\delta\geq2$. If
$G$ contains a connected positive factor $F$ with $d_F(v)<d_G(v)$ for every vertex $v$ of $G$, then $G$ has the strong parity property.

\medskip

A spanning $k$-tree of a connected graph $G$ is a spanning tree in which each vertex admits degree at most $k$. Fan et al. \cite{FSHL} gave a spectral
radius condition for the existence of a spanning $k$-tree in a connected graph.

\medskip

\noindent{\textbf{Theorem 1.3}} (Fan et al. \cite{FSHL}). Let $k\geq3$ be an integer, and let $G$ be a connected graph of order $n\geq2k+16$. If
$\rho(G)\geq\rho(K_1\vee(K_{n-k-1}\cup kK_1))$, then $G$ contains a spanning $k$-tree unless $G=K_1\vee(K_{n-k-1}\cup kK_1)$.

\medskip

A combination of Theorem 1.2 and Theorem 1.3 gives a sufficient spectral condition for a graph to have the strong parity property.

\medskip

\noindent{\textbf{Corollary 1.4.}} Let $k\geq3$ be an integer, and let $G$ be a connected graph of order $n$ with minimum degree $\delta\geq2$.
If $n\geq2k+16$ and
$$
\rho(G)\geq\rho(K_1\vee(K_{n-k-1}\cup kK_1)),
$$
then $G$ has the strong parity property unless $G=K_1\vee(K_{n-k-1}\cup kK_1)$.

\medskip

For $\delta\geq3$, if $k\geq\delta(\delta-2)+1$, then $K_1\vee(K_{n-k-1}\cup kK_1)$ is a spanning subgraph of $H(n,\delta)$, and so
$\rho(K_1\vee(K_{n-k-1}\cup kK_1))<\rho(H(n,\delta))$. By a direct computation, we obtain $\rho(K_1\vee(K_{n-4}\cup 3K_1))>\rho(H(n,\delta))$ for
$\delta\geq3$. Thus, it is easy to see that there exists a positive integer $k_0$ such that
$\rho(K_1\vee(K_{n-k-1}\cup kK_1))>\rho(H(n,\delta))$ for $3\leq k<k_0$ and $\rho(K_1\vee(K_{n-k-1}\cup kK_1))\leq\rho(H(n,\delta))$ for $k\geq k_0$. Furthermore,
$2k+16\geq2\delta^2$ when $k$ is sufficiently large. Consequently, Theorem 1.1 and Corollary 1.4 cannot be obtained from each other. When $k<k_0$, our spectral
radius condition is better than that of Corollary 1.4; when $k$ is sufficiently large, our order is better than that of Corollary 1.4.

\section{Some preliminaries}

Lu, Yang, and Zhang \cite{LYZ} presented a characterization for a graph to have the strong parity property, which is shown in the following.

\medskip

\noindent{\textbf{Lemma 2.1}} (Lu, Yang, and Zhang \cite{LYZ}). A graph $G$ has the strong parity property if and only if
$$
c(G-S)\leq\sum\limits_{v\in S}{d_G(v)}-2|S|+1,
$$
for any $S\subseteq V(G)$.

\medskip

The following lemma follows directly from [2, Theorem 6.8].

\medskip

\noindent{\textbf{Lemma 2.2}} (Bapat \cite{B}). Let $G$ be a connected graph, and let $H$ be a proper subgraph of $G$. Then $\rho(G)>\rho(H)$.

\medskip

\noindent{\textbf{Lemma 2.3}} (Fan, Lin, and Lu \cite{FLLt}, Fan et al. \cite{FSHL}). Let $n_1\geq n_2\geq\cdots\geq n_t\geq p$ be $(t+1)$ positive
integers with $n=\sum_{i=1}^{t}n_i+s$ and $n_1<n-s-p(t-1)$. Then
$$
\rho(K_s\vee(K_{n_1}\cup K_{n_2}\cup\cdots\cup K_{n_t}))<\rho(K_s\vee(K_{n-s-p(t-1)}\cup(t-1)K_p)).
$$

\medskip

Let $M$ be a real matrix whose rows and columns are indexed by $V=\{1,2,\ldots,n\}$. Assume that $M$, with respect to the partition
$\pi:V=V_1\cup V_2\cup\cdots\cup V_r$, can be written as
\begin{align*}
M=\left(
  \begin{array}{cccc}
    M_{11} & M_{12} & \cdots & M_{1r}\\
    M_{21} & M_{22} & \cdots & M_{2r}\\
    \vdots & \vdots & \ddots & \vdots\\
    M_{r1} & M_{r2} & \cdots & M_{rr}\\
  \end{array}
\right),
\end{align*}
where $M_{ij}$ denotes the submatrix (block) of $M$ obtained by rows in $V_i$ and columns in $V_j$. Let $b_{ij}$ denote the average row sum of
$M_{ij}$, that is, $b_{ij}$ is the sum of all entries in $M_{ij}$ divided by the number of rows. Then matrix $M_{\pi}=(b_{ij})$ is called the
quotient matrix of $M$. If the row sum of every block $M_{ij}$ is a constant, then the partition is equitable.

\medskip

\noindent{\textbf{Lemma 2.4}} (Brouwer, and Haemers \cite{BH}, You et al.\cite{YYSX}). Let $M$ be a real matrix with an equitable
partition $\pi$, and let $M_{\pi}$ be the corresponding quotient matrix. Then the eigenvalues of $M_{\pi}$ are also eigenvalues of $M$. Furthermore,
if $M$ is nonnegative and irreducible, then the spectral radius of $M_{\pi}$ is equal to the spectral radius of $M$.

\medskip

The subsequent lemma is the well-known Cauchy Interlacing Theorem.

\medskip

\noindent{\textbf{Lemma 2.5}} (Haemers \cite{Hi}). Let $M$ be a Hermitian matrix of order $s$, and let $N$ be a principal submatrix of $M$ with
order $t$. If $\lambda_1\geq\lambda_2\geq\cdots\geq\lambda_s$ are the eigenvalues of $M$ and $\mu_1\geq\mu_2\geq\cdots\geq\mu_t$ are the eigenvalues
of $N$, then $\lambda_i\geq\mu_i\geq\lambda_{s-t+i}$ for $1\leq i\leq t$.

\section{Proof of Theorem 1.1}

In this section, we first present two lemmas which will be used to prove Theorem 1.1. For integers $\delta, s, n$ with
$4\leq\delta+1\leq s\leq\frac{n-2}{\delta-1}$ and $n\geq2\delta^{2}$, we have the following result on the spectral radii of $H(n,\delta,s)$ and $H(n,\delta)$.

\medskip

\noindent{\textbf{Lemma 3.1.}} For integers $\delta, s, n$ with $4\leq\delta+1\leq s\leq\frac{n-2}{\delta-1}$ and $n\geq2\delta^{2}$, we have
$$
\rho(H(n,\delta,s))<\rho(H(n,\delta)).
$$

\medskip

\noindent{\it Proof.} For the graph $H(n,\delta,s)$, the quotient matrix of $A(H(n,\delta,s))$ corresponding to the partition
$V(H(n,\delta,s))=V(K_s)\cup V(K_{n-(\delta-1)s-1})\cup V(((\delta-2)s+1)K_1)$ is given as
\begin{align*}
B_1=\left(
  \begin{array}{ccc}
    s-1 & n-(\delta-1)s-1 & (\delta-2)s+1\\
    s & n-(\delta-1)s-2 & 0\\
    s & 0 & 0\\
  \end{array}
\right),
\end{align*}
whose characteristic polynomial is
\begin{align*}
\varphi_{B_1}(x)=&x^{3}+(-n+(\delta-2)s+3)x^{2}+(-n-(\delta-2)s^{2}+(\delta-3)s+2)x\\
&+(\delta-2)s^{2}n+sn-(\delta-1)(\delta-2)s^{3}-(3\delta-5)s^{2}-2s.
\end{align*}
Note that the partition $V(H(n,\delta,s))=V(K_s)\cup V(K_{n-(\delta-1)s-1})\cup V(((\delta-2)s+1)K_1)$ is equitable. By virtue of Lemma 2.4, the
largest root, say $\theta_1$, of $\varphi_{B_1}(x)=0$ equals $\rho(H(n,\delta,s))$. Let $\theta_1=\rho(H(n,\delta,s))\geq\theta_2\geq\theta_3$ be
the three roots of $\varphi_{B_1}(x)=0$ and $Q_1=\mbox{diag}(s,n-(\delta-1)s-1,(\delta-2)s+1)$. It is easy to check that
\begin{align*}
Q_1^{\frac{1}{2}}B_1Q_1^{-\frac{1}{2}}=\left(
  \begin{array}{ccc}
    s-1 & s^{\frac{1}{2}}(n-(\delta-1)s-1)^{\frac{1}{2}} & s^{\frac{1}{2}}((\delta-2)s+1)^{\frac{1}{2}}\\
    s^{\frac{1}{2}}(n-(\delta-1)s-1)^{\frac{1}{2}} & n-(\delta-1)s-2 & 0\\
    s^{\frac{1}{2}}((\delta-2)s+1)^{\frac{1}{2}} & 0 & 0\\
  \end{array}
\right)
\end{align*}
is symmetric, and also contains
\begin{align*}
\left(
  \begin{array}{ccc}
    n-(\delta-1)s-2 & 0\\
    0 & 0\\
  \end{array}
\right)
\end{align*}
as its submatrix. Since $Q_1^{\frac{1}{2}}B_1Q_1^{-\frac{1}{2}}$ and $B_1$ have the same eigenvalues, according to Lemma 2.5, we conclude
\begin{align}\label{eq:3.1}
\theta_2\leq n-(\delta-1)s-2\leq n-(\delta-1)(\delta+1)-2=n-\delta^{2}-1.
\end{align}

For the graph $H(n,\delta)$, its adjacency matrix $A(H(n,\delta))$ has the quotient matrix $B_*$ which is obtained by replacing $s$ with $\delta$
in $B_1$, and $B_*$ has the characteristic polynomial $\varphi_{B_*}(x)$ which is obtained by replacing $s$ with $\delta$ in $\varphi_{B_1}(x)$.
Thus, we get
\begin{align}\label{eq:3.2}
\varphi_{B_*}(x)=&x^{3}+(-n+\delta(\delta-2)+3)x^{2}+(-n-\delta^{2}(\delta-2)+\delta(\delta-3)+2)x\nonumber\\
&+(\delta-2)\delta^{2}n+\delta n-\delta^{3}(\delta-1)(\delta-2)-\delta^{2}(3\delta-5)-2\delta.
\end{align}
According to Lemma 2.4, the largest root, say $\eta_*$, of $\varphi_{B_*}(x)=0$ equals $\rho(H(n,\delta))$. Notice that $K_{n-\delta(\delta-2)-1}$
is a proper subgraph of $H(n,\delta)$, it follows from \eqref{eq:3.1}, $\delta\geq3$ and Lemma 2.2 that
\begin{align}\label{eq:3.3}
\eta_*=\rho(H(n,\delta))>\rho(K_{n-\delta(\delta-2)-1})=n-\delta(\delta-2)-2>n-\delta^{2}-1\geq\theta_2.
\end{align}
We are to verify $\varphi_{B_1}(\eta_*)=\varphi_{B_1}(\eta_*)-\varphi_{B_*}(\eta_*)>0$. By a simple computation, we obtain
\begin{align}\label{eq:3.4}
\varphi_{B_1}(\eta_*)=\varphi_{B_1}(\eta_*)-\varphi_{B_*}(\eta_*)=(s-\delta)g(\eta_*),
\end{align}
where $g(\eta_*)=(\delta-2)\eta_*^{2}+(-(\delta-2)s-\delta^{2}+3\delta-3)\eta_*+(\delta-2)(s+\delta)n+n-(\delta-1)(\delta-2)(s^{2}+\delta s+\delta^{2})
-(3\delta-5)(s+\delta)-2$. Then the symmetry axis of $g(\eta_*)$ is $\eta_*=\frac{(\delta-2)s+\delta^{2}-3\delta+3}{2(\delta-2)}$, which implies
that $g(\eta_*)$ is increasing in the interval $\Big[\frac{(\delta-2)s+\delta^{2}-3\delta+3}{2(\delta-2)},+\infty\Big)$. By \eqref{eq:3.3}, $s\geq\delta+1$
and $n\geq(\delta-1)s+2$, we easily prove that
$$
\frac{(\delta-2)s+\delta^{2}-3\delta+3}{2(\delta-2)}<n-\delta(\delta-2)-2<\eta_*,
$$
and so
\begin{align}\label{eq:3.5}
g(\eta_*)>&g(n-\delta(\delta-2)-2)\nonumber\\
=&(\delta-2)n^{2}+(\delta-2)(-2\delta^{2}+4\delta-3)n-(\delta-1)(\delta-2)s^{2}\nonumber\\
&+(-\delta^{2}+\delta+1)s+\delta^{5}-6\delta^{4}+14\delta^{3}-18\delta^{2}+13\delta-4\nonumber\\
\geq&(\delta-2)n^{2}+(\delta-2)(-2\delta^{2}+4\delta-3)n\nonumber\\
&-(\delta-1)(\delta-2)\left(\frac{n-2}{\delta-1}\right)^{2}+(-\delta^{2}+\delta+1)\left(\frac{n-2}{\delta-1}\right)\nonumber\\
&+\delta^{5}-6\delta^{4}+14\delta^{3}-18\delta^{2}+13\delta-4 \ \ \ \ \ \ \ \ \ \ \ (\mbox{since} \ s\leq\frac{n-2}{\delta-1})\nonumber\\
=&\frac{1}{\delta-1}((\delta-2)^{2}n^{2}+(-2\delta^{4}+10\delta^{3}-20\delta^{2}+22\delta-13)n+\delta^{6}\nonumber\\
&-7\delta^{5}+20\delta^{4}-32\delta^{3}+33\delta^{2}-23\delta+10)\nonumber\\
\geq&\frac{1}{\delta-1}(4\delta^{4}(\delta-2)^{2}+2\delta^{2}(-2\delta^{4}+10\delta^{3}-20\delta^{2}+22\delta-13)+\delta^{6}\nonumber\\
&-7\delta^{5}+20\delta^{4}-32\delta^{3}+33\delta^{2}-23\delta+10) \ \ \ \ \ \ \ \ \ \ \ (\mbox{since} \ n\geq2\delta^{2})\nonumber\\
=&\frac{1}{\delta-1}(\delta^{6}-3\delta^{5}-4\delta^{4}+12\delta^{3}+7\delta^{2}-23\delta+10)\nonumber\\
>&0 \ \ \ \ \ \ \ \ \ \ \ (\mbox{since} \ \delta\geq3).
\end{align}

In terms of \eqref{eq:3.4}, \eqref{eq:3.5} and $s\geq\delta+1$, we conclude
$$
\varphi_{B_1}(\eta_*)=(s-\delta)g(\eta_*)>0.
$$
As $\theta_2\leq n-\delta^{2}-1<\rho(H(n,\delta))=\eta_*$ (see \eqref{eq:3.3}), we conclude that $\rho(H(n,\delta,s))=\theta_1<\eta_*=\rho(H(n,\delta))$.
Lemma 3.1 is proved. \hfill $\Box$

\medskip

For positive integers $\delta, s, n$ with $s\leq\delta-1$, $n\geq2\delta^{2}$ and $n\geq((\delta-2)s+2)(\delta+1-s)+s$, we have the following result
on the spectral radii of $H_s$ and $H(n,\delta)$, where $H_s=K_s\vee(K_{n-s-((\delta-2)s+1)(\delta+1-s)}\cup((\delta-2)s+1)K_{\delta+1-s})$.

\medskip

\noindent{\textbf{Lemma 3.2.}} For positive integers $\delta, s, n$ with $\delta\geq3$, $s\leq\delta-1$, $n\geq2\delta^{2}$ and
$n\geq((\delta-2)s+2)(\delta+1-s)+s$, we have
$$
\rho(H_s)<\rho(H(n,\delta)).
$$

\medskip

\noindent{\it Proof.} For the graph $H_s$, consider the partition $V(H_s)=V(K_s)\cup V(K_{n-s-((\delta-2)s+1)(\delta+1-s)})\cup V(((\delta-2)s+1)K_{\delta+1-s})$.
The corresponding quotient matrix of $A(H_s)$ is equal to
\begin{align*}
B_2=\left(
  \begin{array}{ccc}
    s-1 & n-s-((\delta-2)s+1)(\delta+1-s) & ((\delta-2)s+1)(\delta+1-s)\\
    s & n-s-((\delta-2)s+1)(\delta+1-s)-1 & 0\\
    s & 0 & \delta-s\\
  \end{array}
\right).
\end{align*}
Then the characteristic polynomial of $B_2$ is
\begin{align*}
\varphi_{B_2}(x)=&x^{3}+((-\delta+2)s^{2}+(\delta^{2}-\delta-2)s-n+3)x^{2}\\
&+((\delta^{2}-3\delta+2)s^{2}-(n+\delta^{3}-2\delta^{2}-2\delta+1)s+\delta n-n-\delta^{2}-2\delta+2)x\\
&+((2-\delta)s^{3}+(\delta^{2}-\delta-3)s^{2}+\delta s+\delta)n-(\delta-2)^{2}s^{5}\\
&+(2\delta^{3}-6\delta^{2}-\delta+10)s^{4}-(\delta^{4}-2\delta^{3}-6\delta^{2}+7\delta+10)s^{3}\\
&-(2\delta^{3}-\delta^{2}-5\delta-5)s^{2}-(\delta^{3}-\delta)s-\delta^{2}-2\delta.
\end{align*}
According to Lemma 2.4, $\rho(H_s)$ is the largest root of $\varphi_{B_2}(x)=0$.

Recall that $\varphi_{B_*}(x)$ (see \eqref{eq:3.2}) is the characteristic polynomial of $H(n,\delta)$ and $\eta_*=\rho(H(n,\delta))$ is the largest root
of $\varphi_{B_*}(x)=0$. By plugging the value $\eta_*$ into $x$ of $\varphi_{B_2}(x)-\varphi_{B_*}(x)$, we obtain
\begin{align}\label{eq:3.6}
\varphi_{B_2}(\eta_*)=\varphi_{B_2}(\eta_*)-\varphi_{B_*}(\eta_*)=(\delta-s)h(\eta_*),
\end{align}
where $h(\eta_*)=(\delta-2)(s-1)\eta_*^{2}+(n-\delta^{2}s+3\delta s-2s+\delta^{2}-4\delta+1)\eta_*+(\delta^{2}-4\delta+4)s^{4}
-(\delta^{3}-2\delta^{2}-5\delta+10)s^{3}+(\delta n-2n-\delta^{2}-3\delta+10)s^{2}+(-\delta n+3n+\delta^{3}-4\delta^{2}+5\delta-5)s
-\delta^{2}n+2\delta n+\delta^{4}-3\delta^{3}+5\delta^{2}-6\delta$. Using \eqref{eq:3.3}, $2\leq s\leq\delta-1$ and $n\geq((\delta-2)s+2)(\delta+1-s)+s$,
we have
$$
-\frac{n-\delta^{2}s+3\delta s-2s+\delta^{2}-4\delta+1}{2(\delta-2)(s-1)}<n-\delta(\delta-2)-2<\eta_*,
$$
which implies that
\begin{align*}
h(\eta_*)>&h(n-\delta(\delta-2)-2)\\
=&((\delta-2)s-\delta+3)n^{2}+((\delta-2)s^{2}-(2\delta^{3}-7\delta^{2}+10\delta-9)s+2\delta^{3}-9\delta^{2}+12\delta-9)n\\
&+(\delta^{2}-4\delta+4)s^{4}-(\delta^{3}-2\delta^{2}-5\delta+10)s^{3}-(\delta^{2}+3\delta-10)s^{2}\\
&+(\delta^{5}-5\delta^{4}+12\delta^{3}-18\delta^{2}+15\delta-9)s-\delta^{5}+6\delta^{4}-13\delta^{3}+18\delta^{2}-16\delta+6.
\end{align*}
Note that
\begin{align*}
-\frac{(\delta-2)s^{2}-(2\delta^{3}-7\delta^{2}+10\delta-9)s+2\delta^{3}-9\delta^{2}+12\delta-9}{2((\delta-2)s-\delta+3)}=&\frac{2\delta^{2}-3\delta+3-s}{2}\\
\leq&\frac{2\delta^{2}-3\delta+1}{2}<2\delta^{2}\leq n,
\end{align*}
and so
\begin{align*}
h(\eta_*)>&h(n-\delta(\delta-2)-2)\\
\geq&((\delta-2)s-\delta+3)(2\delta^{2})^{2}+((\delta-2)s^{2}-(2\delta^{3}-7\delta^{2}+10\delta-9)s+2\delta^{3}-9\delta^{2}+12\delta-9)(2\delta^{2})\\
&+(\delta^{2}-4\delta+4)s^{4}-(\delta^{3}-2\delta^{2}-5\delta+10)s^{3}-(\delta^{2}+3\delta-10)s^{2}\\
&+(\delta^{5}-5\delta^{4}+12\delta^{3}-18\delta^{2}+15\delta-9)s-\delta^{5}+6\delta^{4}-13\delta^{3}+18\delta^{2}-16\delta+6\\
=&\delta^{3}((s-1)\delta^{2}+s\delta-s^{3}+2s^{2}-8s+11)+(s^{4}+2s^{3}-5s^{2})\delta^{2}\\
&+(-4s^{4}+5s^{3}-3s^{2}+15s-16)\delta+4s^{4}-10s^{3}+10s^{2}-9s+6\\
\geq&\delta^{3}((s-1)(s+1)^{2}+s(s+1)-s^{3}+2s^{2}-8s+11)+(s^{4}+2s^{3}-5s^{2})(s+1)^{2}\\
&+(-4s^{4}+5s^{3}-3s^{2}+15s-16)(s+1)\\
&+4s^{4}-10s^{3}+10s^{2}-9s+6 \ \ \ \ \ \ \ \ \ \ \ \ \ \ \ \ \ \ (\mbox{since} \ \delta\geq s+1\geq3)\\
=&\delta^{3}(4s^{2}-8s+10)+s^{6}+5s^{4}-16s^{3}+17s^{2}-10s-10\\
>&0 \ \ \ \ \ \ \ \ \ \ \ \ \ \ \ \ \ \ (\mbox{since} \ \delta\geq s+1\geq3).
\end{align*}
Combining this with \eqref{eq:3.6} and $\delta\geq s+1$, we obtain
\begin{align}\label{eq:3.7}
\varphi_{B_2}(\eta_*)=(\delta-s)h(\eta_*)>0.
\end{align}

By a simple computation of the first derivative, we have
\begin{align*}
\varphi_{B_2}'(x)=&3x^{2}+2((-\delta+2)s^{2}+(\delta^{2}-\delta-2)s-n+3)x\\
&+(\delta^{2}-3\delta+2)s^{2}-(n+\delta^{3}-2\delta^{2}-2\delta+1)s+\delta n-n-\delta^{2}-2\delta+2.
\end{align*}
Notice that $n\geq((\delta-2)s+2)(\delta+1-s)+s$ and
$$
-\frac{(-\delta+2)s^{2}+(\delta^{2}-\delta-2)s-n+3}{3}<n-\delta(\delta-2)-2,
$$
and so $\varphi_{B_2}'(x)$ is increasing in the interval $[n-\delta(\delta-2)-2,+\infty)$. When $x>n-\delta(\delta-2)-2$, we have
\begin{align}\label{eq:3.8}
\varphi_{B_2}'(x)>&\varphi_{B_2}'(n-\delta(\delta-2)-2)\nonumber\\
=&n^{2}+((-2\delta+4)s^{2}+(2\delta^{2}-2\delta-5)s-4\delta^{2}+9\delta-3)n+(2\delta^{3}-7\delta^{2}+9\delta-6)s^{2}\nonumber\\
&+(-2\delta^{4}+5\delta^{3}-2\delta^{2}-2\delta+7)s+3\delta^{4}-12\delta^{3}+17\delta^{2}-14\delta+2\nonumber\\
\geq&(2\delta^{2})^{2}+((-2\delta+4)s^{2}+(2\delta^{2}-2\delta-5)s-4\delta^{2}+9\delta-3)(2\delta^{2})\nonumber\\
&+(2\delta^{3}-7\delta^{2}+9\delta-6)s^{2}+(-2\delta^{4}+5\delta^{3}-2\delta^{2}-2\delta+7)s\nonumber\\
&+3\delta^{4}-12\delta^{3}+17\delta^{2}-14\delta+2 \ \ \ \ \ \ \ \ \ \ \ \ \ \ \ \ \ \ (\mbox{since} \ n\geq2\delta^{2})\nonumber\\
=&\delta^{2}((2s-1)\delta^{2}+(-2s^{2}+s+6)\delta+s^{2}-12s+11)+(9s^{2}-2s-14)\delta-6s^{2}+7s+2\nonumber\\
\geq&\delta^{2}((2s-1)(s+1)^{2}+(-2s^{2}+s+6)(s+1)+s^{2}-12s+11)\nonumber\\
&+(9s^{2}-2s-14)(s+1)-6s^{2}+7s+2 \ \ \ \ \ \ \ \ \ \ \ \ \ \ \ \ \ \ (\mbox{since} \ \delta\geq s+1\geq3)\nonumber\\
=&\delta^{2}(3s^{2}-5s+16)+(9s^{2}-2s-14)(s+1)-6s^{2}+7s+2\nonumber\\
\geq&(s+1)^{2}(3s^{2}-5s+16)+(9s^{2}-2s-14)(s+1)-6s^{2}+7s+2 \ \ \ \ \ (\mbox{since} \ \delta\geq s+1\geq3)\nonumber\\
=&3s^{4}+10s^{3}+10s^{2}+18s+4\nonumber\\
>&0.
\end{align}
For $s=1$, one can check that $h(\eta_*)>0$, $\varphi_{B_2}(\eta_*)=(\delta-s)h(\eta_*)>0$ and $\varphi_{B_2}'(x)>\varphi_{B_2}'(n-\delta(\delta-2)-2)>0$.
Hence, inequalities \eqref{eq:3.7} and \eqref{eq:3.8} also hold for $s=1$. By \eqref{eq:3.7} and \eqref{eq:3.8}, we deduce
$$
\rho(H_s)<\eta_*=\rho(H(n,\delta)).
$$
This completes the proof of Lemma 3.2. \hfill $\Box$

\medskip

In what follows, we are ready to prove Theorem 1.1.

\medskip

\noindent{\it Proof of Theorem 1.1.} Suppose, to the contrary, that a connected graph $G$ does not have the strong parity property. In terms of Lemma 2.1,
we have
$$
c(G-S)\geq\sum\limits_{v\in S}{d_G(v)}-2|S|+2\geq(\delta-2)|S|+2,
$$
for some nonempty proper subset $S$ of $V(G)$. Let $|S|=s$ and $c(G-S)=t$. Then
$$
t\geq(\delta-2)s+2
$$
and $G$ is a spanning subgraph of $G_1=K_s\vee(K_{n_1}\cup K_{n_2}\cup\cdots\cup K_{n_{(\delta-2)s+2}})$, where $n_1,n_2,\ldots,n_{(\delta-2)s+2}$
are positive integers with $n_1\geq n_2\geq \ldots\geq n_{(\delta-2)s+2}$ and $\sum\limits_{i=1}^{(\delta-2)s+2}{n_i}=n-s$. According to Lemma 2.2,
we obtain
\begin{align}\label{eq:3.9}
\rho(G)\leq\rho(G_1),
\end{align}
with equality if and only if $G=G_1$. We are to proceed by considering the following three possible cases.

\noindent{\bf Case 1.} $s\geq\delta+1$.

For $G_1$ and $H(n,\delta,s)$, according to Lemma 2.3, we conclude
\begin{align}\label{eq:3.10}
\rho(G_1)\leq\rho(H(n,\delta,s)),
\end{align}
where the equality holds if and only if $(n_1,n_2,\ldots,n_{(\delta-2)s+2})=(n-(\delta-1)s-1,1,\ldots,1)$.

It follows from \eqref{eq:3.9}, \eqref{eq:3.10}, $n\geq2\delta^{2}$, $n\geq(\delta-1)s+2$ and Lemma 3.1 that
$$
\rho(G)\leq\rho(G_1)\leq\rho(H(n,\delta,s))<\rho(H(n,\delta)),
$$
which is a contradiction to $\rho(G)\geq\rho(H(n,\delta))$.

\noindent{\bf Case 2.} $s=\delta$.

According to Lemma 2.3, we obtain
$$
\rho(G_1)\leq\rho(H(n,\delta)),
$$
with equality if and only if $G_1=H(n,\delta)$. Together with \eqref{eq:3.9}, we
have
$$
\rho(G)\leq\rho(H(n,\delta)),
$$
where the equality holds if and only if $G=H(n,\delta)$, a contradiction.

\noindent{\bf Case 3.} $s\leq\delta-1$.

Recall that $H_s=K_s\vee(K_{n-s-((\delta-2)s+1)(\delta+1-s)}\cup((\delta-2)s+1)K_{\delta+1-s})$. Note that $G$ is a spanning subgraph of
$G_1=K_s\vee(K_{n_1}\cup K_{n_2}\cup\cdots\cup K_{n_{(\delta-2)s+2}})$, where $n_1,n_2,\ldots,n_{(\delta-2)s+2}$
are positive integers with $n_1\geq n_2\geq \ldots\geq n_{(\delta-2)s+2}$ and $\sum\limits_{i=1}^{(\delta-2)s+2}{n_i}=n-s$. Clearly, the
minimum degree of $G_1$ is at least $\delta$, and so $n_{(\delta-2)s+2}\geq\delta+1-s$. By virtue of Lemma 2.3, we get
\begin{align}\label{eq:3.11}
\rho(G_1)\leq\rho(H_s),
\end{align}
with equality if and only if $(n_1,n_2,\ldots,n_{(\delta-2)s+2})=(n-s-((\delta-2)s+1)(\delta+1-s),\delta+1-s,\ldots,\delta+1-s)$.
According to \eqref{eq:3.9}, \eqref{eq:3.11}, $n\geq2\delta^{2}$, $n\geq((\delta-2)s+2)(\delta+1-s)+s$ and Lemma 3.2, we conclude
$$
\rho(G)\leq\rho(G_1)\leq\rho(H_s)<\rho(H(n,\delta)),
$$
which contradicts $\rho(G)\geq\rho(H(n,\delta))$. Theorem 1.1 is verified. \hfill $\Box$

\section*{Data availability statement}

My manuscript has no associated data.

\section*{Declaration of competing interest}

The authors declare that they have no conflicts of interest to this work.

\section*{Acknowledgments}

The authors are very grateful to the anonymous reviewers for their valuable comments and corrections which result in an improvement of the
original manuscript. This work is supported by the Natural Science Foundation of Jiangsu Province, China. Project ZR2023MA078 supported by
Shandong Provincial Natural Science Foundation.

\end{document}